\documentclass[12pt]{article}

\usepackage{amsmath} 
\usepackage{amsthm} 
\usepackage{amsfonts} 
\usepackage{amssymb}

\newtheorem{thm}{Theorem}[section] 
\newtheorem{defi}[thm]{Definition} 
 
\newtheorem{lemma}[thm]{Lemma}

\usepackage[dvips]{graphicx} 
 
\begin{document} 
 
\title{\bf  Isotropic systems and the interlace polynomial}

\author{Joanna A. Ellis-Monaghan$^1$\\{\small Department of Mathematics}\\ 
{\small Saint Michael's College}\\ 
{\small 1 Winooski Park, Colchester, VT 05439}\\ 
{\small Phone:(802)654-2660, Fax:(802)654-2610}\\ 
{\small jellis-monaghan@smcvt.edu}\\ 
\\Irasema Sarmiento\\ 
{\small Dipartimento di Matematica}\\ 
{\small Universit\'a di Roma ``Tor Vergata"}\\ 
{\small Via della Ricerca Scientifica, I-00133, Rome, Italy}\\ 
{\small Phone: (39) 06 7259-4837, Fax:(39) 06 7259-4699}\\ 
{\small sarmient@mat.uniroma2.it}}\date{} \maketitle 
\noindent 
{\bf Proposed running head:} ISOTROPIC SYSTEMS AND INTERLACE POLYNOMIAL\\

$1$ This research supported by the Vermont Genetics Network through the 
NIH Grant Number 
1 P20 RR16462 from the INBRE  program of the NCRR.

\newpage 
 \begin{abstract}

Through a series of papers in the 1980's, Bouchet introduced 
isotropic systems and the 
Tutte-Martin polynomial of an isotropic system.  Then, Arratia, Bollob\'as, and Sorkin developed 
the interlace polynomial of a 
graph in \cite{ABS00} in response to a DNA sequencing application. 
The interlace polynomial has 
generated considerable recent attention, with new results including 
realizing the original interlace polynomial 
by a closed form generating function 
expression instead of by the original recursive definition 
(see  Aigner and van der Holst \cite{AvdH04}, 
and Arratia, Bollob\'as, and Sorkin \cite{ABS04b}).  Now, 
Bouchet \cite{Bou} recognizes the vertex-nullity interlace polynomial 
of a graph as 
the Tutte-Martin polynomial 
of an associated isotropic system.  This suggests that the machinery of isotropic 
systems may be 
well-suited to investigating properties of the interlace polynomial. 
Thus, we present here an alternative proof for the closed 
form presentation of the vertex-nullity interlace polynomial 
using the machinery of isotropic systems. 
This approach both illustrates the intimate 
connection between the vertex-nullity interlace polynomial and the 
Tutte-Martin polynomial of an isotropic 
system and also provides a concrete example of manipulating 
isotropic systems. We also provide a brief survey of related work.\\

 
 
\end{abstract} 
\footnotetext{ISOTROPIC SYSTEMS AND INTERLACE POLYNOMIAL.} \newpage 
 
\section{Introduction}\label{introduction}

\footnotetext{ISOTROPIC SYSTEMS AND INTERLACE POLYNOMIAL.} 
 
In the 1980's Bouchet unified essential properties of 4-regular graphs 
and pairs of dual binary matriods through a new algebraic construct, 
isotropic systems.  With the Tutte-Martin polynomial of an isotropic 
system, he then significantly extended a fundamental relationship 
between the Tutte polynomial of a planar graph and the Martin polynomial 
of its medial graph.  Then, at the turn of the millennium, Arratia, 
Bollobas, and Sorkin developed the interlace polynomial of a graph to 
analyze the interlaced repeated subsequences of nucleotides that can 
complicate DNA sequencing.  Remarkably, despite the very different 
terminologies, motivations and approaches, the original interlace 
polynomial of a graph may be realized as the Tutte-Martin polynomial of 
an associated isotropic system (see Bouchet, \cite{Bou}).   
 
        The intimate connection between the Tutte-Martin polynomial of 
an isotropic system and the interlace polynomial of a graph suggests 
that the machinery of isotropic systems may be well-suited to 
investigating properties of the interlace polynomial, properties of 
particular current interest due to the interlace polynomial's important 
applications in the biological sciences and general string 
reconstruction problems.  In the current paper we present an alternative 
proof for the closed form presentation of the original vertex-nullity 
interlace polynomial by interpreting it as the Tutte-Martin polynomial 
of a carefully chosen graphical isotropic system.  Our approach places 
the emphasis on manipulating the Tutte-Martin polynomial of an isotropic 
system, using it more concretely than typical in prior work, and thus 
providing an explicit and detailed example of how the machinery of 
isotropic systems may be used to extract information about the interlace 
polynomial.  En route, we survey related work, providing a consolidated 
introduction to the interlace and Tutte-Martin polynomials and their 
interconnections.

A series of papers throughout the 1980's and 1990's, 
including work by Bouchet \cite{Bou87a}, 
\cite{Bou87b}, \cite{Bou88}, 
\cite{Bou89}, \cite{Bou91}, 
\cite{Bou93}, as well as Bouchet and Ghier \cite{BG96}, and 
Jackson \cite{Jac91}, 
developed the theory of isotropic 
systems and of Tutte-Martin 
polynomials of isotropic systems, which 
are now emerging as valuable tools for extending matroidal 
and graph theoretical results.  In addition to the above-mentioned 
works, see also in particular 
Aigner \cite{Aig00}, Aigner and Mielke \cite{AM00}, 
Aigner and van der Holst \cite{AvdH04}, 
Allys \cite{All94}, 
as well as Bouchet's series on multimatroids \cite{Bou97}, \cite{Bou98a}, 
\cite{Bou01}, 
\cite{Bou98b}. 
 
In \cite{ABS00}, Arratia, Bollob\'as, and Sorkin defined a one-variable graph polynomial 
$q_N(G)$ (denoted $q$ there, 
but we follow later work \cite{ABS04b}, reserving $q$ for the 
two-variable generalization) motivated by 
questions arising from DNA sequencing by hybridization in 
Arratia, Bollob\'as, Coppersmith, and Sorkin \cite{ABCS00}. 
This polynomial models the 
interlaced repeated subsequences of DNA that can interfere 
with the unique reconstruction of the original 
DNA strand.  This work promptly generated further interest and 
other applications, for example in 
Aigner, and van der Holst \cite{AvdH04}; 
Arratia, Bollob\'as, and Sorkin \cite{ABS04a}, \cite{ABS04b}; 
Balister, Bollob\'as, Cutler, and Peabody \cite{BBCP02}; 
Balister, Bollob\'as, Riordan, and Scott \cite{BBRS01}; 
Parker and Riera \cite{RP}; and \cite{E-MS}.  In 
\cite{ABS04b},  Arratia, Bollob\'as, and Sorkin define a much richer two-variable interlace 
polynomial, and show that the original polynomial of 
\cite{ABS00} and \cite{ABS04a} is a specialization of it, renaming the 
original interlace 
polynomial the vertex-nullity interlace polynomial due 
to its relationship with the two-variable generalization.

\footnotetext{ISOTROPIC SYSTEMS AND INTERLACE POLYNOMIAL.}

The vertex-nullity interlace polynomial was originally defined recursively via a 
`pivot' operation, and considerable 
cleverness 
was 
required to show that it was in fact well-defined, that is, 
independent of the order of the pivot operations 
(in \cite{ABS00} and \cite{ABS04a}, 
Arratia, Bollob\'as, and Sorkin describe using a computer 
search to ferret out the necessary identities).  One 
of the most significant recent results, derived 
from properties of the two-variable interlace polynomial 
by Arratia, Bollob\'as, and Sorkin in 
\cite{ABS04b} and shown with a slightly different formulation by Aigner and van der Holst in \cite{AvdH04}, 
gives the original vertex-nullity interlace 
polynomial as a closed form generating function expression, 
and thus an alternative proof that it is well-defined. 
Both \cite{ABS04b} and \cite{AvdH04} use linear algebra techniques, 
specifically interpreting 
the pivot operation in terms 
of its impact on adjacency matrices.  Here, we derive the same 
result using purely the properties of isotropic systems 
and the Tutte-Martin polynomial of isotropic systems.

Bouchet   \cite{Bou} has recently explicitly described 
how the recursive form of the vertex-nullity interlace polynomial of a graph can be 
derived from the restricted Tutte-Martin polynomial of an 
isotropic system, a connection implicit by combining 
previous work by Bouchet 
\cite{Bou87a}, \cite{Bou91}, \cite{Bou88}. 
We exploit this relationship and 
the machinery of isotropic systems by viewing the 
interlace polynomial as the Tutte-Martin polynomial 
of an isotropic system with a particular graphic presentation. 
With this, we give an alternative proof for the closed 
form presentation of the vertex-nullity interlace polynomial. 
Although our result, that 
$q(G;x) = \sum_{W\subseteq V(G)}(x-1)^{|W|-r(A(W))}$, 
is precisely the formulation given 
by Arratia, Bollob\'as, and Sorkin in \cite{ABS04b} 
rather than that by Aigner and van der Holst 
in \cite{AvdH04},  where the sum 
involves admissible column sets in an extended adjacency matrix, 
the latter's approach is somewhat closer, although the 
reverse, of that presented here.

In \cite{AvdH04}, Aigner and van der Holst use linear algebra techniques 
to obtain a closed form expression for the vertex-nullity 
interlace polynomial and then compare it to a closed form 
presentation of the Tutte-Martin polynomial of an isotropic 
system to conclude that vertex-nullity interlace polynomial 
and Tutte-Martin polynomial coincide.  Here, we reverse 
this approach, using the main result of Bouchet \cite{Bou} to first 
write $q_N (G;x)$ as the Tutte-Martin polynomial of an isotropic system. 
Then, by choosing the particular isotropic system carefully, 
we are able to conclude the result of 
Arratia, Bollob\'as, and Sorkin \cite{ABS04b}, 
that $q(G;x) = \sum_{W\subseteq V(G)}(x-1)^{|W|-r(A(W))}$.

\section{A brief overview}\label{preliminaries} 
 
The vertex-nullity interlace polynomial of a graph was defined 
recursively by  Arratia, Bollob\'as, and Sorkin in \cite{ABS00} 
via a pivoting operation and was seen by them in \cite{ABS04b} to be 
a specialization of a two-variable interlace 
polynomial, $q(G;x,y)$, with a similar 
pivot recursion.  Let $vw$ be an edge of a graph $G$, and let $A_v$, 
$A_w$ and $A_{vw}$ be the sets of vertices of $G$ adjacent to $v$ only, $w$ only, 
and to both $v$ and $w$, respectively.  The pivot operation ``toggles" 
the edges between $A_v$,  $A_w$ and $A_{vw}$, by deleting existing edges and inserting 
edges between previously non-adjacent vertices. 
The result of this operation 
is denoted $G^{vw}$.  More formally,  $G^{vw}$  has the same vertex set as $G$, and edge 
set equal to the symmetric difference $E(G)\Delta   S$, where $S$  is the 
complete tripartite graph with vertex classes 
$A_v$,  $A_w$ and $A_{vw}$. See Figure $1$.\\

\footnotetext{ISOTROPIC SYSTEMS AND INTERLACE POLYNOMIAL.}

[Insert Figure 1: The pivot operation.] 
 
\begin{defi}\label{vertexnullity} 
 
The vertex-nullity interlace polynomial is defined recursively as: 
 
$$ 
q_N(G;x)= 
\begin{cases} 
x^n {\text{ if }}G=E_n, {\text{ the edgeless graph on $n$ vertices}},\\ 
q_N(G- v;x)+q_N(G^{vw}- w;x) {\text{ if }}vw\in E(G). 
\end{cases} 
$$ 
 
\end{defi}

This polynomial was shown to be well-defined on all simple graphs 
by Arratia, Bollob\'as, and Sorkin 
in \cite{ABS00}.

\begin{defi} 
 
(Arratia, Bollob\'as, and Sorkin \cite{ABS04b}) The two-variable interlace polynomial is defined, 
for a graph $G$ of order $n$, 
by 
 
\begin{equation}\label{redq}q(G;x,y)=\sum_{S\subseteq V(G)}(x-1)^{r(G[W])}(y-1)^{n(G[W])},\end{equation} 
 
\noindent 
where $r(G[W])$ and $n(G[W])=|W|-r(G[W])$ are, respectively, 
the $GF(2)$-rank and 
nullity of the adjacency 
matrix of $G[W]$, the subgraph of $G$ induced by $W$. 
 
\end{defi} 
 
Equivalently, the two-variable interlace polynomial can be defined by the 
following reduction formulas from 
Arratia, Bollob\'as, and Sorkin \cite{ABS04b}. 
 
For a graph $G$, for any edge $ab$ where neither $a$ nor $b$ has 
a loop,

\begin{equation}\label{2eq} 
q(G)=q(G- a)+q(G^{ab}- b)+((x-1)^2-1)q(G^{ab}- a- b), \\ 
\end{equation} 
 
\footnotetext{ISOTROPIC SYSTEMS AND INTERLACE POLYNOMIAL.} 
 
for any looped vertex $a$, 
 
$$q(G)=q(G- a)+(x-1)q(G^a- a),$$ 
 
\noindent 
and, for the edgeless graph $E_n$ on $n\geq 0$ vertices, $q(E_n)=y^n$. 
Here $G^a$ is the {\em local complementation} of $G$ , and is defined 
as follows. Let $N(a)$ be the neighbors of $a$, that is, the set 
$\{w\in V: a {\text{ and }}w {\text{ are joined by an edge}}\}$. 
Thus $a\in N (a)$ iff $a$ is a loop. 
The graph $G^a$ is equal to $G$ except 
that $G^a[N (a)] =\overline{G[N (a)]}$, i.e. we ``toggle" the edges among the neighbors of $a$, switching edges to non-edges and vice-versa.

Arratia, Bollob\'as, and Sorkin show in \cite{ABS04b} that the vertex-nullity interlace polynomial is a 
specialization of the two-variable interlace 
polynomial as follows:

\begin{equation}\label{stateqN} 
q_N(G;y)=q(G;2,y)=\sum_{W\subseteq V(G)}(y-1)^{n(G[W])}. 
\end{equation} 
 
An equivalent formulation for $q_N(G;x)$  is given by Aigner and van der Holst in 
\cite{AvdH04}: $q_N(G;g)=\sum_S(y-1)^{co(L_S)}$, where 
the sum is over certain colum sets of $L$, an extended adjacency 
matrix of $G$.

In \cite{ABS00} and \cite{ABS04b}, the Arratia, Bollob\'as, and Sorkin 
give an interpretation 
of the vertex-nullity interlace polynomial of a 
circle graph in terms of the circuit partition, or Martin, polynomial of a related 
$4$-regular Eulerian digraph.  A {\em circle graph} 
on $n$ vertices is a graph $G$ 
derived from a chord diagram, 
where two copies of each of the symbols $1$ through $n$ are 
arranged on the perimeter of a circle, 
and a chord is drawn between like symbols. 
Two vertices 
$v$ and $w$ in $G$ share an edge if and only if their corresponding 
chords intersect in the chord diagram. 
See Figure $2$. 
 
Circle graphs have also been called alternance graphs 
by Bouchet in \cite{Bou88} and 
interlace graphs by Arratia, Bollob\'as, and Sorkin  in \cite{ABS00}. 
Research on circle graphs includes a 
complete characterization and a polynomial time algorithm for 
identifying them. For example, see 
Bouchet \cite{Bou85}, \cite{Bou87b}, \cite{Bou87c}, \cite{Bou94}; 
Czemerinski, Duran, Gravano, and Bouchet \cite{CDG02}; 
Dur\'an \cite{Dur03}; 
de Fraysseix \cite{Fra84}; Gasse   \cite{Gas97}; 
Read and  Rosenstiehl  \cite{RR78a}, \cite{RR78b}; and 
Wessel and  P\"oschel \cite{WP84}.\\

%\noindent 

%\\ 
 
[Insert Figure 2: The circle graph of a chord diagram.]\\

A $4$-{\em regular Eulerian digraph} is a $4$-regular directed graph such that, 
at each vertex, two edges are oriented inward, and two are oriented outward. 
A $4$-regular Eulerian digraph is called a $2$-in, $2$-out graph in 
Arratia, Bollob\'as, and  Sorkin \cite{ABS00}. Note that if 
$C$ is an Eulerian circuit of a $4$-regular Eulerian digraph, and we write the vertices 
along the perimeter of a circle in the order that they are visited by $C$ (each is visited 
exactly $2$ times), and then draw a chord between like vertices, the result is a chord diagram.

\footnotetext{ISOTROPIC SYSTEMS AND INTERLACE POLYNOMIAL.}

\begin{defi} 
\rm 
A {\em graph state} of a $4$-regular Eulerian digraph $\vec{G}$  is the result of replacing 
each $4$-valent vertex $v$ of  $\vec{G}$ with two $2$-valent vertices each 
joining an incoming and an outgoing 
edge originally adjacent to $v$.  Thus a graph state is a disjoint union of consistently 
oriented cycles.  See Figure $3$. 
\end{defi}

Note that graph states (see \cite{E-M98}) are equivalent to the circuit 
partitions of Arratia, Bollob\'as, and  Sorkin  \cite{ABS00} and 
Bollob\'as \cite{Bol02}, the Eulerian 
decompositions of Bouchet \cite{Bou88}, and the Eulerian 
$k$-partitions of Martin \cite{Mar77} and Las Vergnas \cite{Las83}.\\

[Insert Figure $3$: A graph state.]\\ 
 
\begin{defi} 
\rm 
The {\em circuit partition polynomial} of a $4$-regular Eulerian 
digraph  $\bar{G}$ is 
$f(\vec{G};x)=\sum_{k\geq 0}f_k(\vec{G})x^k$, where $f_k(\vec{G})$ 
is the number of graph states of $\bar{G}$  with $k$ components, defining 
$f_0(\vec{G})$ to be $1$ if $\vec{G}$  has no edges, 
and $0$ otherwise. 
\end{defi} 
 
The circuit partition polynomial is a simple translation 
of the Martin polynomial $m(\vec{G};x)$, 
defined 
recursively for $4$-regular digraphs by Martin in his 1977 thesis 
\cite{Mar77}, with 
$f(\vec{G};x)=xm(\vec{G};x+1)$.

Martin also showed that if 
$G$ is a planar graph and $\vec{G_m}$ is an appropriately oriented 
medial graph of $G$, then $m(\vec{G_m};x)=T(G;x,x)$, where 
$T(G;x,y)$ is the Tutte polynomial, giving a catalyst for the 
subsequent development of the Tutte-Martin polynomial. 
 
Las Vergnas found closed forms for the Martin polynomials (for both 
graphs and digraphs). He also extended their properties to general 
Eulerian digraphs and further developed their theory 
(see \cite{Las79}, \cite{Las88}, \cite{Las83}). The transforms of the 
Martin polynomials, for arbitrary Eulerian graphs and digraphs, 
were given in \cite{E-M98}, and then aptly named circuit partition 
polynomials by Bollob\'as in \cite{Bol02}, with splitting 
identities provided in \cite{Bol02} and \cite{E-Mb}. The circuit partition 
polynomial is also a specialization of a much broader 
multivariable polynomial, the generalized transition polynomial of 
\cite{E-MS02}, which assimilates such graph invariants as the 
Penrose polynomial that are not evaluations of the 
Tutte polynomial.

\footnotetext{ISOTROPIC SYSTEMS AND INTERLACE POLYNOMIAL.}

For circle graphs, the vertex-nullity 
interlace polynomial and the circuit partition polynomial 
are related 
by the following theorem.

\begin{thm}(Arratia, Bollob\'as, and  Sorkin \cite{ABS00}, Theorem 6.1).\label{A} 
 
If $\vec{G}$  is a $4$-regular Eulerian digraph, $C$ is any Eulerian 
circuit of $\vec{G}$ , and  $H$ is the circle graph of the chord diagram 
determined by $C$, 
 then $f(\vec{G};x)=xq_N(H;x+1)$ . 
 
\end{thm}

We note that this result is also indicated (in terms of the original Martin 
polynomial not the circuit 
partition polynomial) using isotropic systems in the work of Bouchet, although the 
components appear in three separate papers.  The necessary results for expressing 
the Martin polynomial of a graph as a Tutte-Martin polynomial of an isotropic system are given 
in \cite{Bou87a}, and an outline for assembling them to recover the Martin polynomial of a 
$4$-regular 
Eulerian (di)graph from the appropriate isotropic system appears in \cite{Bou91}. 
Then, in \cite{Bou88}, it is proved that every graphic isotropic system (i.e. those associated 
in a special way to $4$-regular graphs) has a circle graph as a fundamental graph (see \cite{Bou88}, 
theorem 6.3 in particular). 
Thus the Tutte-Martin polynomial of the graphic isotropic system 
that is equal to the Martin polynomial of a $4$-regular Eulerian digraph 
$\vec{G}$ is in turn equal 
to the Tutte-Martin polynomial of the isotropic system of the circle graph of an Eulerian 
circuit of $\vec{G}$  . 
 
The following brief introduction to isotropic systems 
follows closely 
that provided by Bouchet in \cite{Bou}. Further details can be found 
for example in 
\cite{Bou87a} and \cite{Bou91}.

Let $K=\{0,x,y,z\}$ be the Klein group, so in particular 
$K$ is isomorphic to ${\mathbb Z}_2\times {\mathbb Z}_2$ with 
$0\leftrightarrow (0,0)$, 
$x\leftrightarrow (1,0)$, 
$y\leftrightarrow (0,1)$, and 
$z\leftrightarrow (1,1)$. We view $K$ as a vector space over 
$GF(2)$, 
the two element field $\{0,1\}$.

There is a bilinear form $<.,.>$ uniquely defined on $K$, mapping 
$K$ into $GF(2)$ with $<a,b>=1$ if $a\neq b$ and neither $a$ nor $b$ 
is zero, and $<a,b>=0$ otherwise. For every finite set $V$, the set 
$K^V$ of maps from $V$ to $K$ is a vector space of dimension $|V|$ 
over $GF(2)$. If the elements of $V$ are ordered, then these maps 
may be thought of as 
vectors of length $|V|$ where the $i^{th}$ entry is value in 
$K$ assigned to the $i^{th}$ element of $V$. 
 
We equip $K^V$ with a bilinear form $<.,.>$ defined by 
$<A,B>=\sum_{v\in V}<A(v),B(v)>$ for all $A$, $B$ in 
$K^V$. 
 
We set $K'=K\setminus \{0\}=\{x,y,z\}$.

\footnotetext{ISOTROPIC SYSTEMS AND INTERLACE POLYNOMIAL.}

\begin{defi} 
An isotropic system is a pair $S=(L,V)$, where $V$ is a finite set and $L$ 
is a vector subspace of $K^V$ such that $dim L= | V |$ and $<A,B>=0$, for 
all $A,B\in L$. 
\end{defi} 
 
Let $G$ be a simple graph with vertex set $V$ and the edge set $E$. The set 
${\mathcal P}(V)$ of subsets of $V$ is considered as a 
vector space over $GF(2)$. 
Thus, if $V$ is ordered, a subset of $V$ may be identified with 
a vector of length $|V|$, with $i^{th}$ entry $1$ if the $i^{th}$ element of 
$V$ is in the subset and $0$ else. Addition of two vectors in the 
vector space thus corresponds to the symmetric difference of the underlying 
subsets.

For $v\in V$, let $N(v)$ be the neighbors of $v$, that is, the set 
$\{w\in V: v {\text{ and }}w {\text{ are joined by an edge}}\}$. 
We set $N(P)=\sum_{v\in P}N(v)$, for $P\subseteq V$. 
Note that $N(P)$, since we are summing in $GF(2)$, 
is simply the set of vertices of $G$ that are adjacent 
to an odd number of the vertices of $P$. 
For $X\in K^V$ and $P\subseteq V$, the vector $Y\in K^V$ defined by 
$Y(v)=X(v)$ if $v\in P$ and $Y(v)=0$ if $v\notin P$, will be 
denoted by $X(P)$. 
 
\begin{thm}(Bouchet \cite{Bou88} Theorem 3.1 )\label{2.5} 
Let $G$ be a simple graph with vertex set $V$. 
Let $A,B\in K'^V$ be such that $A(v)\neq B(v)$, for all $v\in V$, 
and set $L=\{A(P)+ B(N(P)):P\subseteq V\}$. Then $S=(L,V)$ is an 
isotropic system.

\end{thm} 
 
The triple $(G,A,B)$ is called a {\em graphic presentation} of $S$. 
 
For $X\in K^V$ we set $\widehat{X}=\{Y\in K^V:Y(v)\in\{0,X(v)\},v\in V\}$. 
Equivalently, $\widehat{X}=\{X(P)\mid P\subseteq V\}$. 
We note that $\widehat{X}$ is a vector subspace of $K^V$ with $X(P)+X(Q)=X(P\Delta Q)$. 
 
\begin{defi}\label{2.7} 
For $C\in K'^V$, the restricted Tutte-Martin polynomial 
$m(S,C;\xi)$ is defined by 
 
$$m(S,C;\xi)=\sum(\xi-1)^{dim(L\cap \widehat{F})},$$

where the sum is over $F\in K'^V$ such that $F(v)\neq C(v)$, 
for all $v\in V.$ 
\end{defi}

\footnotetext{ISOTROPIC SYSTEMS AND INTERLACE POLYNOMIAL.} 
 
Defintion \ref{2.7} gives the Tutte-Martin polynomial of a general 
isotropic system. However, in the special case that the 
isotropic system is given by a graphic presentation, we have 
the following central result.

\begin{thm}(Bouchet \cite{Bou} Equality (5))\label{Bou2.8} 
If $G$ is a simple graph and $S$ is the isotropic system defined 
by a graphic 
presentation $(G,A,B)$, then 
 
$$q_N(G;\xi)=m(S,A+B;\xi).$$ 
 
\end{thm} 
 
At the heart of the proof is that $q_N(G;\xi)$ may be recursively defined by 
the following relations: 
 
\begin{itemize} 
 
\item[(1)]  $q_N(G;\xi)=1$ if $G$ is the empty graph. 
 
\item[(2)]  $q_N(G;\xi)=\xi q_N(G\setminus v;\xi)$ if $v$ is an isolated 
vertex,

\item[(3)]  $q_N(G;\xi)= q_N(G\setminus v;\xi)+q_N(G*vwv\setminus v;\xi)$ 
if $vw\in E$.

\end{itemize}

Here $G*v$ is the local complementation of $G$ at $v$, which replaces 
the subgraph induced by $N(v)$ by the complementary subgraph, precisely as defined above, but the 
notation of Bouchet \cite{Bou} differs from 
Arratia, Bollob\'as, and  Sorkin \cite{ABS04b}. 
$G*vwv$ is the iterated operation $G*v*w*v$. See Figure 
$4$, and compare to Figure $1$.\\ 
 
[Insert Figure $4$: complementation along an edge].\\

\footnotetext{ISOTROPIC SYSTEMS AND INTERLACE POLYNOMIAL.}

Aigner and van der Holst  give an alternative proof of this same 
result in \cite{AvdH04}. If $G$ is a simple graph on $n$ 
vertices, then let $L$ be the $(n\times 2n)$ matrix over 
$GF(2)$ whose first $n$ columns are the adjacency matrix of $G$, and 
whose last $n$ columns are the $n\times n$ identity matrix, labeling 
the rows of $L$ by $1,\dots ,n$, and the columns by 
$1,\dots, n, \bar{1},\dots ,\bar{n}$. A column set $S$ is 
{\em admissible} if $|S\cap \{i,\bar{i}\}|=1$ for all $i$, and $L_S$ is the 
$(n\times n)$ submatrix with column set $S$. Having first used 
a linear algebra approach to show that 
$q_N(G;x)=\sum_S (x-1)^{co(L_S)}$, where the sum is over all 
admissible column sets $S$, 
and $co(L_S)$ is the corank, Aigner and van 
der Holst then compare this to the 
formula of Definition \ref{2.7} to conclude 
Theorem \ref{Bou2.8}. 
 
\section{Interpretation of $q_N(G)$ via Isotropic systems}

That $q_N(G;x) = m(S, C; x)$ gives an interpretation of 
what $q_N(G;x)$ counts, but only with respect to the underlying isotropic system. 
Here, by translating the isotropic machinery into graph theoretical terms, 
we give an alternative proof of the generating function formulation 
of Arratia, Bollob\'as, and  Sorkin \cite{ABS04b}, namely that 
 
\begin{equation}\label{2} 
q_N(G;x) = \sum_{W\subseteq V(G)}(x-1)^{\mid W\mid-r(M(W))}, 
\end{equation}

\noindent 
where $M(W)$ 
is the adjacency matrix of $G$ restricted to $W$ 
and $r$ is the $GF(2)$ rank function.

Given a simple graph $G$ with vertex set $V$, we choose an isotropic system 
with a particular graphic presentation as follows.  For $p\in \{x,y,z\}$ 
we define $\bar{p}:V\rightarrow K$ as the element of $K^V$ such that $\bar{p}(v)=p$ 
for all $v\in V$.  Let $A=\bar{x}$, $B=\bar{y}$, and 
$C=A+B=\bar{z}$. Now let $S$ be the isotropic system whose graphic presentation is 
$(G,A,B)$ as given by Theorem \ref{2.5}. 
 
\footnotetext{ISOTROPIC SYSTEMS AND INTERLACE POLYNOMIAL.} 
 
Using Definition \ref{2.7} and recalling that 
and $\widehat{F}=\{F'\in K^V \mid F'(v)=F(v) {\text{ or }}F'(v)=0\}$, we 
have, for $A=\bar{x}$, $B=\bar{y}$, that 
 
$$m(S,\bar{z};\xi)=\sum(\xi-1)^{dim(L\cap\widehat{W})},$$ 
 
where 
the sum is over $F\in\{x,y\}^V$, 
and where $L=\{A(P)+B(N(P)):P\subseteq V\}$. 
 
Write $L_P$ for $A(P)+B(N(P))$, so, 
 
\begin{equation}\label{eq5} 
L_P(v)= 
\begin{cases} 
x {\text{ if }} v\in P, v\notin N(P),\\ 
y {\text{ if }} v\notin P, v\in N(P),\\ 
z {\text{ if }} v\in P\cap N(P),\\ 
0 {\text{ else. }} 
\end{cases} 
\end{equation}

If $L_P\in L\cap\widehat{F}$ and $F\in \{x,y\}^V$, then $L_P\in\{0,x,y\}^V$. 
For $P\subseteq V$, we denote by $G\mid P$ 
the subgraph of $G$ induced by $P$, and we write $F_x$ for $F^{-1}(x)$ and $F_y$ for $F^{-1}(y)$. An {\em even} graph has all vertices of 
even degree and is not necessarily connected.

\begin{lemma}\label{A1} 
For $P\subseteq V$ the following statements are equivalent: 
\begin{enumerate} 
\item 
$L_P\in\{0,x,y\}^V$, 
 
\item 
$G\mid P$ is an even subgraph of $G$. 
 
\end{enumerate} 
\end{lemma} 
\begin{proof} 
It follows from $(\ref{eq5})$ that 
$L_P\in\{0,x,y\}^V$ if and only if $P\cap N(P)=\emptyset$. 
Since $N(P)$ is the set of vertices adjacent to an odd number of 
vertices of $P$, this is equivalent to saying that every 
$v\in P$ has an even number of neighbors in $P$. That is, if 
and only if $G\mid P$ is an even subgraph of $G$. 
 
\end{proof}

Since Aigner and van der Holst 
\cite{AvdH04} and the current work, although giving reverse approaches, are 
in fact examining the same phenomenon, it is not surprising that 
they should at some point converge, and we note that the preceding 
Lemma \ref{A1} is essentially Lemma $3$ of \cite{AvdH04}.

\begin{lemma}\label{B1} 
Let $F\in\{x,y\}^V$ and $P\subseteq V$. Then the following statements are equivalent: 
\begin{enumerate} 
\item 
$L_P\in \widehat{F}$, 
 
\item 
$G\mid P$ is an even subgraph of $G\mid F_x$ and $N(P)\subseteq F_y$. 
 
\end{enumerate}

\end{lemma} 
\begin{proof} 
 
By the definition of $L_P$ we have that $L_P\in \widehat{F}$ if and only if there is 
and 
$R\subseteq V$ such that 
 
\begin{equation}\label{3} 
L_P(v)=F(R)= 
\begin{cases} 
x {\text{ if }} v\in R\cap F_x\\ 
y {\text{ if }} v\in R\cap F_y\\ 
0 {\text{ if }} v\notin R 
\end{cases} 
\end{equation} 
 
Comparing (\ref{3}) to (\ref{eq5}), we note that $L_P=F(R)$ if and only if $P=R\cap F_x$, $N(P)=R\cap F_y$ 
and 
$R=P\cup N(P)$. That is, $L_P=F(R)$ if and only if $P\subseteq F_x$, $N(P)\subseteq F_y$, 
$P\cap N(P)=\emptyset$ and $R=P\cup N(P)$. 
 
By Lemma \ref{A1}, $P\cap N(P)=\emptyset$ and $F\subseteq F_x$ if and only if 
$G\mid P$ is an even subgraph of $G\mid F_x$. 
 
\footnotetext{ISOTROPIC SYSTEMS AND INTERLACE POLYNOMIAL.} 
 
The desired conclusion follows readily. 
 
\end{proof}

By Lemma \ref{B1}, we have that 
 
$$L\cap \widehat{F}=\{L_P\in L: G\mid P {\text{ is an even subgraph of }}G\mid F_x 
{\text{ and }}N(P)\subseteq F_y\}.$$ 
 
Next we determine, in terms of the graph $G$, the dimension of $L\cap\widehat{F}$ as a subspace of 
$K^V$. 
 
Write $star(v)$ for the subgraph of edges incident to $v$. 
Note that $N(P)\subseteq F_y=V\setminus F_x$ if and only if none of the elements 
of $F_x\setminus P$ is in $N(P)$. That is, $N(P)\subseteq F_y$ if and only if for all $v\in F_x\setminus P$ 
the cardinality of $star(v)\cap P$ is even. On the other hand, $G\mid P$ is an even subgraph of $G\mid F_x$ if and only if 
$P\subseteq F_x$ and for all $v\in P$ the cardinality of $star(v)\cap P$ is even. Therefore, by Lemma 
\ref{B1},

\footnotetext{ISOTROPIC SYSTEMS AND INTERLACE POLYNOMIAL.}

$$L\cap \widehat{F}=\{L_P:P\subseteq F_x {\text{ and for all }}v\in F, {\text{ the cardinality of }}star(v)\cap P 
{\text{ is even}}\}.$$ 
 
Recall that the power set ${\mathcal P}(V)$ is a vector space over $GF(2)$ where addition is the symmetric 
differences of subsets of $V$. Let $U$ be the subspace 
of ${\mathcal P}(V)$ defined by 
$U=\{P\subseteq F_x: {\text{ for all }}v\in F, {\text{ the cardinality of }}star(v)\cap P 
{\text{ is even }}\}.$

\begin{lemma}\label{4.3} 
Let $F$ be an element of $\{x,y\}^V$. Then $dim(L\cap\widehat{F})=dim(U)$. 
\end{lemma} 
\begin{proof} 
 
Let $\Psi:L\cap \widehat{F}\rightarrow U$ be the function given by $\Psi (L_P)=P$. 
 
It is straightforward to prove that $\Psi$ is a well defined bijection. 
 
We will prove that $\Psi$ is also a homomorphism. Let $P, P'$ be arbitrarily given subsets 
of $F_x$. As in Lemma \ref{B1}, $L_P=F(P\cup N(P))$, where the union is disjoint. Note that, by the 
definition of the neighborhood of a set, 
$N(P)+N(P')= \sum_{v\in P}N(v)+\sum_{v\in P'}N(v)=\sum_{v\in P+P'}N(v)=N(P+P')$. 
 
Since $L\cap\widehat{F}$ is a subspace of $\widehat{F}$, and we noted below Theorem 
\ref{2.5} that $F(R)+F(Q)=F(R\Delta Q)=F(R+Q)$, we have the following: 
 
$$L_P+L_{P'}=F(N\cup N(P))+F(P'\cup N(P´))$$ 
 
$$=F(P\cup N(P))\Delta F(P'\cup N(P'))$$ 
 
$$=F(P\Delta P')\cup (N(P)\Delta N(P'))$$ 
 
$$=F((P\Delta P') \cup N(P\Delta P'))$$ 
 
$$=L_{P\Delta P'}=L_{P+P'}.$$ 
 
Thus $\Psi$ is a homomorphism. 
 
It follows that $dim(L\cap\widehat{F})=dim(U)$. 
 
\end{proof}

\begin{lemma}\label{4.4} 
Let $F$ be an element of $\{x,y\}^V$. Then 
 
$$dim(L\cap\widehat{F})=dim(U)=\mid F_x\mid-r(M),$$ 
 
where $M$ is the adjacency matrix of $G\mid F_x$. 
\end{lemma}

\footnotetext{ISOTROPIC SYSTEMS AND INTERLACE POLYNOMIAL.}

\begin{proof} 
Let $\psi:\mathcal{P}(F_x)\rightarrow GF(2)^{F_x}$ be the homomorphism whose matrix $M$ is the 
adjacency matrix of $G\mid F_x$. That is, the rows and columns of $M$ are indexed by the elements 
of $G\mid F_x$ and $M_{ij}=1$ iff $i$ and $j$ are joined by an edge in $G\mid F_x$. Also 
$\psi(P)=M\cdot \chi_P$, where $\chi_P$ is the 
vector identified with 
$P$ written 
as a column. 
 
We have that $P\in Ker(\psi)$ if and only if $v(P)=M\cdot \chi_P=0$ (the all zero 
column). That is, for every $i\in F_x$, $0=\sum_{j\in F_x}M_{ij}\chi_{P_j}=\sum_{j\in P}M_{ij}$. 
In other words, for every $i\in F_x$, the cardinality of $star(v)\cap P$ is even. 
 
Thus $L\cap\widehat{F}$ is the kernel of $\psi$ and $dim(Im(\psi))=r(M)$, 
where $r$ is the $GF(2)$ rank function. 
On the other hand $dim({\mathcal P}(F_x))=\mid F_x\mid.$ 
 
Therefore 
 
$$dim(L\cap\widehat{F})=dim(U)=\mid F_x\mid-r(M).$$ 
 
\end{proof}

\begin{thm}\label{4.5} 
Let $G$ be a simple graph. Let $q_N(G;x)$ be the vertex-nullity interlace polynomial of $G$. 
Then 
 
$$q_N(G;x)=\sum_{W\subseteq V(G)}(x-1)^{\mid W\mid-r(M(W))},$$ 
 
where $M(W)$ is the adjacency matrix of $G\mid W$. 
 
\end{thm} 
\begin{proof} 
By Theorem \ref{Bou2.8}, given a simple graph $G$ on the vertex set $V$, 
supplementary vectors 
$A$ and $B$ of $K^V$, and $L=\{A(N(P))+B(P):P\in {\mathcal P}(V)\}$, 
then the pair $S=(L,V)$ is an isotropic system. 
Moreover, the mapping $X\rightarrow A(N(P))+B(P)$ is a 
linear bijection from ${\mathcal P }(V)$ onto $L$. 
Let $A=\bar{x},B=\bar{y}$ and $C=A+B=\bar{z}$ as before. It follows 
that there is an isotropic system $S$ such that 
$(G,A,B)$ is its graphic presentation. 
 
From Definition \ref{2.7} and Theorem 
\ref{Bou2.8},

$$q_N(G;\xi)=m(S,C;\xi)= \sum(\xi-1)^{dim(L\cap\widehat{F})},$$ 
 
\noindent 
where the sum is over $F\in\{x,y,x\}^V$ such that $F(v)\neq C(v)$ for all $v\in V$.

\footnotetext{ISOTROPIC SYSTEMS AND INTERLACE POLYNOMIAL.}

Therefore, for $A=\bar{x}, B=\bar{y}$ and $C=A+B=\bar{z}$, we get 
 
$$m(S,C;\xi)=\sum (\xi -1)^{dim(L\cap\widehat{F})},$$ 
 
where the sum is over $F\in\{x,y\}^V$. 
 
By Lemma \ref{4.4}, $dim(L\cap\widehat{F})=dim(U)=\mid F_x\mid -r(M(F_x))$, 
where $M(F_x)$ is the adjacency matrix of $G\mid F_x$ and $F_x=F^{-1}(x).$ 
Thus,

$$m(S,C;\xi)=\sum_{F_x\subseteq V} (\xi -1)^{\mid F_x\mid -r(M(F_x))}.$$

Therefore, 
 
$$q_N(G;x)=\sum_{W\subseteq V(G)}(x-1)^{\mid W\mid -r(M(W))},$$ 
 
where $M(W)$ is the adjacency matrix of $G\mid W$.

\end{proof}

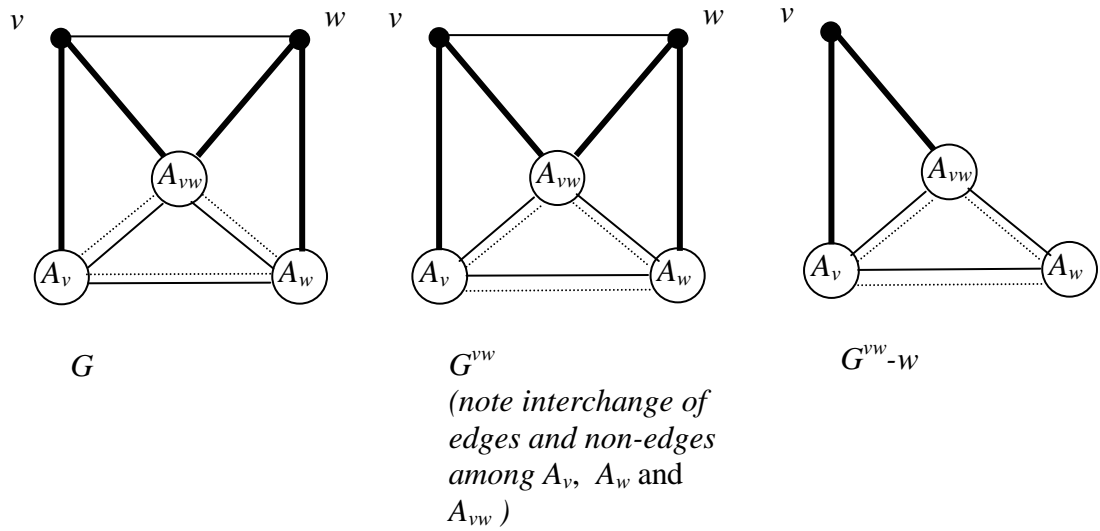
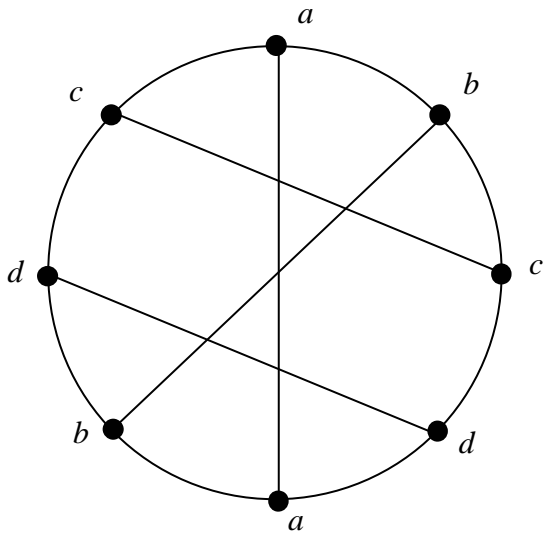
\end{document}